\RequirePackage{fix-cm}
\documentclass[smallextended]{svjour3}       
\smartqed  
\usepackage{graphicx}
\usepackage{amssymb}
\begin{document}

\title{High order structure preserving explicit methods for solving
linear-quadratic optimal control problems and differential games
}


\author{Sergio Blanes
}


\institute{Sergio Blanes \at
              Instituto de Matem\'{a}tica Multidisciplinar, Universitat Polit\`{e}cnica de Val\`{e}ncia.\\
 46022 Valencia, Spain \\
              \email{serblaza@imm.upv.es}           
}

\date{Received: date / Accepted: date}

\maketitle

\begin{abstract}
We present high order explicit geometric integrators to solve
linear-quadratic optimal control problems and $N$-player
differential games. These problems are described by a system
coupled non-linear differential equations with boundary
conditions. We propose first to integrate backward in time the
non-autonomous matrix Riccati differential equations and next to
integrate forward in time the coupled system of equations for the
Riccati and the state vector. This can be achieved by using
appropriate splitting methods, which we show they preserve most
qualitative properties of the exact solution. Since the coupled
system of equations is usually explicitly time dependent, a
preliminary analysis has to be considered. We consider the time as
two new coordinates, and this allows us to integrate the whole
system forward in time using splitting methods while preserving
the most relevant qualitative structure of the exact solution. If
the system is a perturbation of an exactly solvable problem, the
performance of the splitting methods considerably improves. Some
numerical examples are also considered which show the performance
of the proposed methods.
\end{abstract}

\section{Introduction}
\label{intro}

Linear-quadratic (LQ) optimal control problems appear in many
different fields in engineering
\cite{abou03mre,anderson90ocl,reid,speyer10poo} as well as in
quantum mechanics \cite{palao02qcb,peirce88oco,zhu98arm} (see also
\cite{brif10coq} and references therein).
In general, LQ optimal control problems are described by coupled
systems of nonlinear differential equations with boundary
conditions. The particular algebraic structure of the equations
makes that, in general, the solutions have some qualitative
properties which are relevant for the theoretical study, and for
this reason we consider its preservation by the numerical methods
to be of great interest in order to get reliable and accurate
results.

Numerical methods for solving non linear BVPs are usually more
involved (and computationally more costly) than initial value
problems  IVPs \cite{ascher88nso,keller76nso,na79cmi}.
 We consider the numerical integration of LQ optimal control
problems with high order explicit and structure preserving
methods, i.e. methods which preserve most qualitative properties
of the problem. These methods are referred as structure preserving
methods. The method are explicit and the algorithms can be used
with variable step as well as variable order, and they have a
simple way to estimate the accuracy of the method.
%

The main idea is that, once the problem is reformulated as an IVP,
to split the vector field into parts which are exactly solvable
and such that each part preserves the desired structure of the
solution. Then, high order splitting methods can be used. If the
problem is explicitly time dependent or it can be considered as a
small perturbation of an exactly solvable problem, the correct
split and the splitting methods to be used require a more careful
analysis. These methods can also be considered as exponential
methods which have shown a high performance for linear problems
\cite{blanes09tme,blanes02hoo,iser,iserles99ots}, and are usually
referred as geometric numerical integrators.

The methods proposed are also adjusted for solving
linear-quadratic $N$-player differential games (which have been
extensively studied from the theoretical point of view
\cite{anderson90ocl,ba95,Cr71,engw,SH69}) since they can be
considered as optimal linear control problems.

The paper is organized as follows. In Section~\ref{LQcontrol} we
introduce the equations to be solved in a LQ optimal control
problem and present the non-linear matrix Riccati equation as an
equivalent linear system of equations. In
Section~\ref{NumericalMethods} we formulate the problem as an IVP
with a previous backward time integration of the Riccati equation,
we present a brief introduction to splitting methods and we show
how these methods can be used on autonomous and non-autonomous
problems. The case of perturbed systems is and the preservation of
some qualitative properties by splitting methods is also
considered. Section~\ref{DifferentialGames} considers the
generalization to the numerical integration of differential games,
and Section~\ref{NumericalExamples} is devoted to numerical
experiments to illustrate the performance of the methods. Finally,
Section~\ref{Conclusions} gives the conclusions of the work.

\section{LQ optimal control problems}
 \label{LQcontrol}

Let us consider the LQ optimal control problem
\begin{equation} \label{LinearControl}
 x' =  A (t) x + B_{} (t) u_{}(t) ,  \quad \ \ x(0) =
 x_0 \, , \quad \ \ 0=t_{0} \leq t \leq t_{f}=T \, ,
\end{equation}
where the unknown $ x(t)\in \mathbb{R}^{n} $ is the dynamic state.
Here
 $ A (t) \in \mathbb{R}^{n \times n} \, $,
 $ B_{} (t) \in \mathbb{R}^{n \times r_{}} $,  and
 $ \ u_{}(t) \in \mathbb{R}^{r_{}}$ is the control.

We consider a quadratic cost function given by
\begin{equation}\label{}
  J_{} = x^{T}(T) \, Q_{T} \, x(T)  \label{cost}
 + \int_{0}^{T} \left\{ x^{T}(t) Q_{} (t) x(t) + u_{}^{T}(t)
R_{} (t) u_{}(t) \right\} dt \, , \nonumber
\end{equation}
where $ Q (t), Q_{T} \in \mathbb{R}^{n \times n}$,  are symmetric
non-negative matrices, $ R_{} (t) \in \mathbb{R}^{r_{} \times
r_{}}$ is symmetric and positive definite (i.e. $ Q (t), Q_{T}\geq
0, \ R>0$) and $z^{T}$ denotes the transpose of $z$. It is well
known that the optimal control is reached when $u$ is written as
\cite{abou03mre,anderson90ocl,speyer10poo}
\begin{equation}\label{controls}
u_{} (t) \, = \, - \, R_{}^{-1} (t) \, B_{}^{T} (t)  \, P_{} (t)
\,  x(t) \, ,
\end{equation}
with $ \, P_{} (t)\in \mathbb{R}^{n \times n} \, $ satisfying the
matrix Riccati differential equation (RDEs) with final conditions
\begin{equation}\label{CoupleRic}
  P' = - Q (t) - A^{T}(t) P - P A (t) + P S(t) P  , \qquad
  P_{}(T)=Q_{T},
\end{equation}
wherein
\begin{equation}\label{S}
S_{} (t) \, = \, B_{} (t) \, R_{}^{-1}  (t) \, B_{}^T  (t) \, ,
\end{equation}
is a symmetric $n \times n $ matrix and $S(t)\geq 0$. It is known
that $P(t)$ is also a symmetric and non-negative matrix.
Substituting (\ref{controls}) into (\ref{LinearControl}) we have
that
\begin{equation} \label{LinearOptControl}
 x' =  \Big( A (t) - S(t) P(t)\Big) x ,  \quad \ \ x(0) =
 x_0 \, .
\end{equation}
Notice that the product $S(t) P(t)$ is a non-negative matrix and
the preservation of this property by numerical methods is
important to stabilize the evolution of the state vector.

The goal is to compute $u(t)$ to be used to control the system,
and this requires to numerically solve the coupled system of non
linear differential equations with boundary conditions
(\ref{CoupleRic}) and (\ref{LinearOptControl}). The numerical
solution of non linear BVPs is usually very involved, with
computationally costly methods and significant storage
requirements (one has to numerically integrate backward in time
some of the equations, to store intermediate results and  to make
a forward integration using the stored results, or to use shooting
methods \cite{ascher88nso,na79cmi}).

%

We present a new algorithm which solves the problem as an IVP and
then it allows to use variable step and variable order methods
while preserving the geometric structure of the problem. In
addition, the methods need low storage requirements and allows an
immediate evaluation of the controls along the integration, being
of great interest for real time control problems.


To this purpose it is useful to write the non-linear matrix RDE as
an equivalent linear differential equation.

\subsection{The matrix Riccati differential
equation}

For solving the RDE (\ref{CoupleRic}), it is usual to consider the
following decomposition $PU=V$, with $U(t),
V_{}(t)\in\mathbb{R}^{n \times n}$. Let us denote
\begin{displaymath}
 y(t)=\left[
 \begin{array}{c} U \\ V_{}
 \end{array} \right];
 \quad K(t)=\left[
 \begin{array}{ccc} A(t) & - S_{}(t) \\
 - Q_{} (t)& - A^{T} (t)  \end{array} \right],
\end{displaymath}
where  $y(t)\in\mathbb{R}^{2n \times n}, \ K(t)\in\mathbb{R}^{2n
\times 2n}$, $S(t)$ is given by (\ref{S}), and the matrices $A(t)$
and $Q_{}(t)$, are given by (\ref{LinearControl}) and
(\ref{cost}), respectively.  Then, it is easy to check that $y(t)$
is the solution of the IVP
  \begin{equation}\label{y}
 y'(t) = K(t)\, y(t) ; \qquad \quad
  y(T)=\left[ \begin{array}{c} U(T) \\ V_{} (T) \end{array} \right]
  = \left[ \begin{array}{c} I \\ Q_{T} \end{array} \right] \, ,
\end{equation}
with conditions at the end of the interval, to be integrated
backward in time. By \cite{engw,JP95}, if (\ref{y}) has an
appropriate solution with $U(t)$ non singular, the solution of
(\ref{CoupleRic}) can be calculated by
\begin{equation}\label{pi}
P_{} (t) \, = \, V(t) \, U(t)^{-1} \,   .
\end{equation}
Conditions under which $U^{-1}$ exists are known (see
\cite{abou03mre,JoPo} and references therein). In this work we
assume $U(t)$ is non-singular, otherwise $P(t)$ would be unbounded
and the equation for the state vector would not be well defined.


\section{Geometric integrators for solving LQ optimal control problems}
 \label{NumericalMethods}

The coupled system of equations to be solve is given by

\begin{eqnarray}
 \frac{d}{dt} \left[
 \begin{array}{c} U \\ V_{}
 \end{array} \right] & = &
 \left[
 \begin{array}{ccc} A(t) & - S_{}(t) \\
 - Q_{} (t)& - A^{T} (t)  \end{array} \right]
 \left[
 \begin{array}{c} U \\ V_{}
 \end{array} \right], \quad  \quad
 \left[ \begin{array}{c} U(T) \\ V_{} (T) \end{array} \right]
  = \left[ \begin{array}{c} I \\ Q_{T} \end{array} \right]
    \label{eq:RiccatiBack}  \\
  \frac{dx}{dt}  & = &  \Big( A(t) - S(t)  \, V(t) \, U(t)^{-1} \,\Big) x ,
  \qquad \ \ x(0) = x_0 \, .  \label{eq:StateForward1}
\end{eqnarray}

Since the equation (\ref{eq:RiccatiBack}) is independent of
(\ref{eq:StateForward1}), we propose to integrate backward the RDE
to get $U(t_0),V(t_0)$ with sufficiently high accuracy. The method
for this backward integration can take large steps without
frequent outputs and no storage requirement are necessary.
Different methods can be used, and the best choice can depend on
each particular problem (i.e. if the equation is autonomous, or if
it is non-autonomous but the time-dependent matrices have a smooth
time variation, etc.).

Intermediate solutions are not needed in this preliminary
integration so, we have much freedom on the choice of the
numerical method for the backward integration. This allows a fast
way to get the initial conditions, $U(t_0),V(t_0)$, for the RDE.
Next, we have to integrate forward in time the following IVP
\begin{eqnarray*}
 \frac{d}{dt} \left[
 \begin{array}{c} U \\ V_{}
 \end{array} \right] & = &
 \left[
 \begin{array}{ccc} A(t) & - S_{}(t) \\
 - Q_{} (t)& - A^{T} (t)  \end{array} \right]
 \left[
 \begin{array}{c} U \\ V_{}
 \end{array} \right], \quad  \quad
 \left[ \begin{array}{c} U(t_0) \\ V_{} (t_0) \end{array} \right]
  = \left[ \begin{array}{c} U_0 \\ V_{0} \end{array} \right]
    \\
  \frac{dx}{dt}  & = &  \Big( A(t) + S(t)  \, V(t) \, U(t)^{-1} \,\Big) x ,
  \qquad \ \ x(0) = x_0 \,   
    \\
 u (t) & = & - \, R^{-1}(t) \, B^{T}(t) \, V(t) \, U(t)^{-1}\, x(t) \, .
\end{eqnarray*}

The linearized RDE and the equation for the state vector can be
written in short as follows
\begin{equation}\label{eq:Separable}
\begin{array}{rcl}
  v' & = & M(t) v \\
  x' & = & N(t,\hat v)x
\end{array}
\end{equation}
with $v=[U^T \ V^T]^T, \ \hat v= VU^{-1}=P$ and $M,N$ are matrices
of appropriate dimensions. This system of equations is separable
into insolvable parts and then splitting methods can be used in a
simple way to solve this problem. These methods preserve most
qualitative properties of the solutions. Let us briefly introduce
the idea of splitting methods as well as the particular methods
which will be used.

\subsection{Splitting Methods}

 Let us consider the initial value problem
\begin{equation}   \label{eq.1.1}
   x' = f(x), \qquad x_0 = x(0) \in \mathbb{R}^D
\end{equation}
with $f: \mathbb{R}^D \longrightarrow  \mathbb{R}^D$ and solution
$x(t)=\varphi_t(x_0)$, and suppose that the vector field is
separable
\begin{equation}
  \label{eq:f=fa+fb}
f(x)=f^{[a]}(x)+f^{[b]}(x),
\end{equation}
in such a way that the equations
\begin{eqnarray}   \label{eq.1.2}
   x' & = & f^{[a]}(x), \qquad x(0) = x_0^a, \\
   x' & = & f^{[b]}(x), \qquad x(0) = x_0^b,
\end{eqnarray}
can be integrated exactly, with solutions $x(h) =
\varphi_h^{[a]}(x_0^a)$ and  $x(h) = \varphi_h^{[b]}(x_0^b)$,
respectively, at $t = h$, the time step. It is well known that the
composition, $\chi_h = \varphi^{[b]}_{h}\circ \varphi^{[a]}_{h}$,
is a first order method and
\begin{equation}\label{eq:leapfrog}
\psi_h =   \varphi^{[a]}_{h/2}\circ \varphi^{[b]}_{h}\circ
 \varphi^{[a]}_{h/2} ,
\end{equation}
is a symmetric second order method. Higher order methods can be
obtained by composition
\begin{equation}  \label{eq:splitting}
   \psi_h = 
   \varphi^{[b]}_{b_{m}h} \circ \varphi^{[a]}_{a_{m}h}\circ \cdots\circ
  \varphi^{[a]}_{a_{2}h}\circ \varphi^{[b]}_{b_{1}h}\circ
 \varphi^{[a]}_{a_{1}h} ,
\end{equation}
for appropriate choices of coefficients $a_i,b_i$. For simplicity,
we denote the composition as : $b_{m}\,a_m\,\ldots\,b_1\,a_1$. For
example, an efficient fourth-order method is given by the
following sequence \cite{blanes02psp}
\[ 
  b_7\,a_7\,b_6\,a_6\,b_5\,a_5\,b_4\,a_4\,b_3\,a_3\,b_2\,a_2\,b_1\,a_1
\]  
 with $a_1=0$ and the scheme is symmetric ($a_{8-i}=a_{i+1}, \ b_{8-i}=b_i, \ i=1,2,\ldots$)
so, it can be written as follows
\begin{equation}\label{4thOrder}
  b_1\,a_2\,b_2\,a_3\,b_3\,a_4\,b_4\,a_4\,b_3\,a_3\,b_2\,a_2\,b_1
\end{equation}
which is a 6-stage method. It has 6 coefficients $b_i$ and 7
coefficients $a_i$, but the last evaluation of one step can be
reused in the first map in the following step, and it is not
counted for the computational cost (it is called the First Same As
Last (FSAL) property). A sixth-order method is given by the
following 10-stage symmetric sequence
\begin{equation}\label{6thOrder}
  a_1\,b_1\,a_2\,b_2\,a_3\,b_3\,a_4\,b_4\,a_5\,b_5\,
  a_6\,b_5\,a_4\,b_4\,a_4\,b_3\,a_3\,b_2\,a_2\,b_1\,a_1.
\end{equation}
The coefficients of both methods (taken from \cite{blanes02psp})
are collected in Table~\ref{tab.1} for convenience of the reader.

\begin{table}[tbp]
\caption{Coefficients for different splitting methods}
\label{tab.1} {
\begin{tabular}{l}
\begin{tabular}{lll}
\hline\hline
\multicolumn{3}{c}{Coefficients of the 6-stage 4th-order method (\ref{4thOrder})} \\
\hline
   $b_1= 0.0792036964311957 $ & $a_2= 0.209515106613362$\\
   $b_2= 0.353172906049774  $ & $a_3=-0.143851773179818$\\
   $b_3=-0.0420650803577195 $ & $a_4= 1/2 - (a_2+a_3)$ \\
   $b_4= 1 - 2(b_1+b_2+b_3)$  &
\end{tabular}
\\
\begin{tabular}{lll}
\hline%
\multicolumn{3}{c}{Coefficients of the 10-stage 6th-order method (\ref{6thOrder})} \\
\hline
$a_1= 0.0502627644003922  $ &$b_1= 0.148816447901042$ \\
$a_2= 0.413514300428344   $ &$b_2=-0.132385865767784$\\
$a_3= 0.0450798897943977  $ &$b_3= 0.067307604692185$\\
$a_4=-0.188054853819569   $ &$b_4= 0.432666402578175  $    \\
$a_5= 0.541960678450780   $ &$b_5= 1/2 - (b_1+\ldots+b_4)$ \\
$a_6= 1 - 2(a_1+\ldots+a_5)$   & \\\hline
\multicolumn{3}{c}{Coefficients of the 2-stage (4,2) method (\ref{42Order})} \\
\hline
$\displaystyle a_1= \frac{3-\sqrt{3}}{6}  $ &$b_1=1/2 $ \\
$a_2= 1 - 2a_1$  \\\hline
\multicolumn{3}{c}{Coefficients of the 5-stage (8,4) method (\ref{84Order})} \\
\hline
$a_1= 0.07534696026989288842  $ &$b_1=0.19022593937367661925 $ \\
$a_2= 0.5179168546882567823  $ &$b_2= 0.84652407044352625706$\\
$a_3= 1/2 - (a_1+a_2)  $ &$b_3=1 - 2(b_1+b_2)$ \\\hline
\end{tabular}
\end{tabular}
}
\end{table}

The splitting methods we have presented are valid for autonomous
equations because the exact solution of the equation
(\ref{eq.1.1}) can be written formally as the exponential a Lie
operator, $\varphi_t=\exp\big(t\ L_f  \big)$, with $L_f\equiv
f\cdot \nabla$. The map $\varphi_t$ can be approximated by a
composition of exponentials of Lie operators associated to the
vector fields $f^{[a]}$ and $f^{[b]}$. If the problem is
non-autonomous, as usually is the case for LQ optimal control
problems, one can take the time as a new coordinate and transform
the original non-autonomous equation into an autonomous one in an
extended phase space. However, for our problem it is advantageous
to consider the time not as one but as two independent
coordinates, and this has to be combined with an appropriate
splitting method. For this reason, we consider separately the
autonomous from the non-autonomous  LQ optimal control problem.

\subsection{The autonomous case}

If the system is autonomous, we can compute the solution of the
RDE (\ref{eq:RiccatiBack}) at $t_0$ as follows.

\begin{equation}\label{eq:SolBackAutonomous}
 \left[
 \begin{array}{c} U_0 \\ V_{0}
 \end{array} \right]  =
 \exp  \left((t_0-T) \left[
 \begin{array}{ccc} A & - S_{} \\
 - Q& - A^{T}  \end{array} \right] \right)
 \left[
 \begin{array}{c} U_T \\ V_{T}
 \end{array} \right],
\end{equation}
which can be accurately computed, for example, by a scaling and
squaring method \cite{almohy10cta,moler03ndw} i.e. if we take
$K_i$ as an accurate approximation to the scaled exponential
\begin{equation}\label{}
 K_i \simeq \exp \left(\frac{(t_0-T)}{2^i} \left[
 \begin{array}{ccc} A & - S \\
 - Q& - A^{T}  \end{array} \right] \right)
\end{equation}
(here $K_i$ can be, for example, a diagonal Pad\'e approximation
or a Taylor approximation) then
\begin{equation}\label{}
 \left[
 \begin{array}{c} U_0 \\ V_{0}
 \end{array} \right]  \simeq
 \underbrace{\left( \cdots \left( K_i  \right)^2 \cdots \right)^2}_{i-times}
 \left[
 \begin{array}{c} U_T \\ V_{T}
 \end{array} \right],
\end{equation}
which provides an accurate approximation to $U_0,V_0$ from
$U_T,V_T$ with no intermediate results. Next, one has to integrate
forward in time the system (\ref{eq:Separable}) which takes the
form:
\begin{equation}\label{eq:SepAut1}
  \frac{d}{dt} \left[
 \begin{array}{c} v \\ x
 \end{array} \right] =
  \left[
 \begin{array}{c} Mv \\  N(\hat v)x
 \end{array} \right]
\end{equation}
which is separable into two solvable parts
\begin{equation}\label{eq:SepAut2}
  \frac{d}{dt} \left[
 \begin{array}{c} v \\ x
 \end{array} \right] =
  \left[
 \begin{array}{c} Mv \\ 0
 \end{array} \right],
 \qquad \qquad
 \frac{d}{dt} \left[
 \begin{array}{c} v \\ x
 \end{array} \right] =
 \left[
 \begin{array}{c} 0 \\ N(\hat v)x
 \end{array} \right].
\end{equation}
The following algorithm can be used to advance one time step, $h$,
from $t_n$ to $t_{n+1}=t_n+h$, by using an splitting method with
coefficients $\{ a_i,b_i\}_{i=1}^{m}$
\[
\begin{array}{l}
  (v^0,x^0)=(v_n,x_n) \\
  {\bf do} \ \ i=1,\ldots,m \\
     \quad {x}^i  = \exp\left(a_ihN(\hat v^{i-1})\right)x^{i-1} \\
     \quad {v}^i  = \exp\left(b_ihM\right)v^{i-1} \\
  {\bf enddo} \\
  (v_{n+1},x_{n+1})=(v^m,x^m)\\
  u_{n+1} = - \, R^{-1} \, B^{T} \, V_{n+1} \, U_{n+1}^{-1}\, x_{n+1} \, .
\end{array}
\]
In order to save computational cost, we can compute and store the
following exponentials
\[
  E_i = \exp\left(b_ihM\right), \qquad i=1,2,\ldots,m.
\]
To take advantage of the FSAL property to save the computation of
a map,  one has to slightly adjust the previous algorithm. Notice
that this algorithm allows to compute the controls, $u_n$, jointly
with the numerical solution of the equations and this can be
convenient for real time experiments. Obviously, at $t=T$ we must
reach the exact solution $P(T)=Q_T$ with accuracy up to round off
accuracy.

\subsection{The non-autonomous case}

If the system is non-autonomous, the solution of the RDE has no
analytic solution in a closed form. In this case we have to solve
numerically the RDE backward in time using, for example, high
order Magnus integrators (e.g. a method of order six or eight
\cite{blanes02hoo}) which have shown a high performance while
preserving most qualitative properties. Other methods like a high
order extrapolation method can also be used. Once we have
$U_0,V_0$, we need to integrate forward in time the system of non
autonomous equations. If we split the system as in the autonomous
case, we end up with two non-autonomous problems, both with no
solution in a closed form. We are then looking for an alternative
split which provides the same accuracy while preserving the
qualitative structure of the problem. This can be achieved  by
frozen the time following a proper sequence. To this purpose, we
take the time not as a new coordinate, but as two new coordinates
as follows

\begin{equation}\label{eq:Separable2t}
  \left\{
\begin{array} {rcl}
 \displaystyle  v' & = & M(t_1) v \\
 \displaystyle  x' & = & N(t_2,\hat v)x \\
 \displaystyle  t_1' & = & 1 \\
 \displaystyle  t_2' & = & 1
\end{array}
  \right.
\end{equation}
where $ \, '\equiv  \frac{d}{dt}, \, $ and which we split as
follows
\begin{equation}\label{eq:Separable2tb}
  \left\{
\begin{array} {rcl}
 \displaystyle  v' & = & M(t_1) v \\
 \displaystyle  x'& = & 0 \\
 \displaystyle  t_1' & = & 0 \\
 \displaystyle  t_2' & = & 1
\end{array}
  \right.
  \qquad \qquad
  \left\{
\begin{array} {rcl}
 \displaystyle  v' & = & 0 \\
 \displaystyle  x' & = & N(t_2,\hat v)x \\
 \displaystyle  t_1' & = & 1 \\
 \displaystyle  t_2' & = & 0
\end{array}
  \right.
\end{equation}
This corresponds to two linear autonomous equations in the
extended phase space, and each part is now exactly solvable. Then,
the same splitting methods as in the autonomous case can be used.
The algorithm for one time step, $h$, is given by:
\begin{equation}\label{eq:algorithm1}
\begin{array}{l}
  (v^0,x^0,t_1^0,t_2^0)=(v_n,x_n,t_{n},t_{n}) \\
  {\bf do} \ \ i=1,m \\
     \quad {x}^i  = \exp\big(\, a_ihN(t_1^{i-1},\hat v^{i-1})\,\big)x^{i-1} \\
     \quad t_2^i  = t_2^{i-1} + a_ih   \\
     \quad {v}^i  = \exp\big(\,b_ihM(t_2^i)\,\big) v^{i-1} \\ 
     \quad t_1^i  = t_1^{i-1} + b_ih   \\
  {\bf enddo} \\
  (v_{n+1},x_{n+1},t_{n+1},t_{n+1})=(v^m,x^m,t_1^m,t_2^m)\\
  u_{n+1} = - \, R^{-1}(t_{n+1}) \, B^{T}(t_{n+1}) \, V_{n+1} \, U_{n+1}^{-1}\, x_{n+1} \, .
\end{array}
\end{equation}
Since the matrix $M(t)$ changes at each step, the exponentials,
$\exp\big(b_ihM(t_2^i)\big)$, need to be computed at each step. If
this is the most time consuming part of the algorithm, it is
possible to look for an algorithm in which the exponentials can be
replaced by a much more economical symmetric second order
approximation, as for example a second order diagonal Pad\'e
approximation
\[
 \displaystyle   \Phi_{b_ih}^{[b]}(t_2^i) = \frac{I+\frac12 b_ihM(t_2^i)}{I-\frac12 b_ihM(t_2^i)}
  =  \exp(b_ihM(t_2^i)) + \mathcal{O}(h^3).
\]
In general, $M$ belongs to the Lie algebra of symplectic matrices
and the exponential belongs to the associated symplectic Lie
group. Diagonal Pad\'e approximations preserve this group property
for the symplectic algebra \cite{iser}. Now it is important to
keep in mind that one can not use the coefficients from
Table~\ref{tab.1} since they are designed for problems which are
separable in two parts such that both are exactly solved.

We can proceed as follows. Let $S^{[2]}_h$ be the following
symmetric second order method
\begin{equation}\label{eq:S2}
 S^{[2]}_h=  \varphi^{[a]}_{h/2} \circ \Phi_{h}^{[b]}\circ \varphi^{[a]}_{h/2} :
 \left\{ \begin{array}{l}
     \quad x^{1/2}  = \exp\left(\frac{h}{2}N(t_1^{0},\hat v^{0})\right)x^{0} \\
     \quad t_2^{1/2}  = t_2^{0} + \frac{h}{2}   \\
     \quad v^1  = \Phi^{[b]}_{h}(t_2^{1/2})v^{i-1} \\
     \quad t_1^1  = t_1^{0} + h   \\
     \quad x^{1}  = \exp\left(\frac{h}{2}N(t_1^{1},\hat v^{1})\right)x^{1/2} \\
     \quad t_2^{1}  = t_2^{1/2} + \frac{h}{2}   \\
\end{array}  \right.
\end{equation}
Taking $S^{[2]}_h$ as the basic method, we can build methods of
order, $p$, with $p>2$, as a composition of this basic method
\begin{equation}\label{eq:ProdS2}
   S^{[p]}_h = \prod_{i=1}^m   S^{[2]}_{\alpha_i h}.
\end{equation}
The following fourth-order method shows a good performance to get
solutions up to relatively high accuracies (for methods built in
this way)
\[
   S^{[4]}_h =  S^{[2]}_{\alpha_1 h}\circ  S^{[2]}_{\alpha_1 h}\circ
   S^{[2]}_{\alpha_2 h}\circ  S^{[2]}_{\alpha_1 h}\circ  S^{[2]}_{\alpha_1 h}
\]
with $\alpha_1=1/(4-4^{1/3}), \ \alpha_2=-4^{1/3}/(4-4^{1/3})$.
Several sets of coefficients for methods of different orders and
number of stages are collected in
\cite{hairer06gni,mclachlan02sm}.

Notice that at $t=T$ we must reach a numerical approximation close
to the exact solution $P(T)=Q_T$. The methods presented in this
work will approximate this value with accuracy up to the order of
the method, i.e.  $P_{ap}(T)=Q_T+\mathcal{O}(h^{p})$, and this can
be used as a measure of the accuracy of the algorithm.

\subsection{Methods for near-integrable problems}

In some cases, LQ optimal control problems can be formulated as a
small perturbation of an exactly solvable problem (or which can be
easily solved by a numerical method). In this case it can be
convenient to split between the dominant part and the
perturbation, and to use splitting methods designed for separable
problems with this structure. Then, if we have the IVP
\begin{equation}   \label{eq.Pert1}
   x' = f^{[a]}(x) + \varepsilon f^{[b]}(x)
\end{equation}
where $|\varepsilon |\ll 1$ and both equations are either exactly
solvable or can be efficiently solved by a numerical method,
splitting methods tailored for this problem usually have a very
good performance. For example, two highly efficient methods are
given by the following compositions and whose coefficients from
\cite{mclachlan95cmi} are collected in Table~\ref{tab.1}
\begin{equation}\label{42Order}
  a_1\,b_1\,a_2\,b_1\,a_1
\end{equation}
which is referred as a (4,2) method (a second order methods which,
in the limit $ \varepsilon\rightarrow 0$, behaves as a fourth
order method) and
\begin{equation}\label{84Order}
  a_1\,b_1\,a_2\,b_2\,a_3\,b_3\,a_3\,b_2\,a_2\,b_1\,a_1
\end{equation}
which corresponds to a (8,4) method. More elaborated methods with
more stages are given in \cite{blanes13nfo}.

If the problem is non-autonomous
\begin{equation}   \label{eq.Pert2}
   x' = f^{[a]}(x,t) + \varepsilon f^{[b]}(x,t)
\end{equation}
it is still possible to take the time as a new coordinate (only
one new coordinate instead of two contrarily to the previous and
more general case) and the structure of a near integrable problem
remains, in the extended phase space, if we split the system as
follows \cite{blanes10sac}
\begin{equation}   \label{eq.Pert3}
   \frac{d}{dt} \left[ \begin{array}{c} x \\ t_1 \end{array} \right] =
    \left[ \begin{array}{c} f^{[a]}(x,t_1) \\ 1  \end{array} \right] +
    \varepsilon  \left[ \begin{array}{c} f^{[b]}(x,t_1) \\ 0 \end{array}
    \right].
\end{equation}
If the time is considered as two different coordinates as in the
general separable case, this structure of a perturbed problem is
lost. This requires to integrate exactly or up to high accuracy
the non-autonomous equation associated to the dominant part
\begin{equation}   \label{eq.Pert4}
   x' = f^{[a]}(x,t),
\end{equation}
and to solve the perturbed part with the time frozen.

For example, LQ optimal control problems where the matrix $A$ is
constant and $\|A\| \gg\|S\|$ and $\|A\| \gg\|Q\|$ have this
structure. We can split the problem as follows
\begin{equation} \label{eq:PertForwardP}
 x' = f^{[a]}(x,t) \ : \ \left\{ \begin{array}{rcl}
 \left[
 \begin{array}{c} U \\ V_{}
 \end{array} \right]' & = &
 \left[
 \begin{array}{ccc} A & 0 \\
 0& - A^{T}  \end{array} \right]
 \left[
 \begin{array}{c} U \\ V_{}
 \end{array} \right],    \\
  x'  & = &  \Big( A - S(t)  \, V(t) \, U(t)^{-1} \,\Big) x  \,
\end{array}  \right.
\end{equation}
with solution for the matrix RDE: $U=e^{(t-t^*)A}U^*, \
V=e^{-(t-t^*)A^T}V^*$, where $U^*, \, V^*$ are the initial
conditions at each stage at $t=t^*$. The equation for the state
vector is then
\[
  x' = \Big( A + S(t) \, e^{-(t-t^*)A^T} V^* \, U^{*-1}e^{-(t-t^*)A} \,\Big)
  x.
\]
This is a linear equation of the form

\begin{equation}\label{eq.Linear}
  x' =  M(t) x
\end{equation}
which in general has no analytical solution in a closed form, but
can be numerically solved to high accuracy by using, for example,
a Magnus integrator \cite{blanes09tme}. Notice that, for the
numerical integration, the matrix $M(t)$ has to be evaluated on a
number of quadrature nodes, $t_i=t^*+c_ih$, for a time step, $h$.
Then, the exponentials $e^{-c_ihA}$ have to be evaluated only once
at the nodes of the quadrature rule and can be used on each time
step.

Magnus integrators are exponential integrators which preserve most
qualitative properties of the exact solution. The methods are
explicit and usually have similar stability properties than
implicit methods because most implicit methods can be seen as
rational approximations to the exponentials appearing in the
Magnus integrators. Standard Magnus integrators \cite{blanes09tme}
involve commutators of the matrix $M(t)$ evaluated a different
instants. Magnus integrators which do not involve commutators
(commutator-free methods) also exist and in many cases are simpler
to implement. For example, a fourth-order method for the
integration from $t_n$ to $t_{n+1}=t_n+h$ is given by
\[
  x_{n+1}=\exp\left( \frac{h}{12} \left(-M_0+4M_{1/2}+3M_1 \right) \right)
  \exp\left( \frac{h}{12} \left(3M_0+4M_{1/2}-M_1 \right) \right)
  x_n,
\]
where $M_i=M(t_n+ih)$. High order commutator-free methods can be
found in \cite{alvermann11hoc,blanes06fas}.

Next, the coordinate $t_1$ is advanced and for the perturbation we
have to solve the autonomous problem ($ x'=\varepsilon
f^{[b]}(x,t_1)$)
\begin{eqnarray}
 \frac{d}{dt} \left[
 \begin{array}{c} U \\ V_{}
 \end{array} \right] & = &
 \left[
 \begin{array}{ccc} 0 & - S_{}(t_1) \\
 - Q_{} (t_1)& 0  \end{array} \right]
 \left[
 \begin{array}{c} U \\ V_{}
 \end{array} \right],     \label{eq:RiccatiForwardP2}
\end{eqnarray}
where the value of  $t_1$ is frozen. Since $\|Q\|$ and  $\|S\|$
are small, it suffices, for most practical purposes, to
approximate the exact solution by a low order Taylor method.

\subsection{Structure preservation of the splitting methods}

It is well known that given the autonomous matrix RDE
\begin{equation}\label{eq:Preserve1}
  P' = - Q - A^{T} P - P A + P S P  , \qquad
  P(T)=Q_{T},
\end{equation}
where $Q,S,Q_T$ are symmetric non negative matrices, the solution
$P(t)$ is also a symmetric matrix and $P(t)\geq 0, \ t\in[0,T]$.
The same result is also valid for the non-autonomous equation
\begin{equation}\label{eq:Preserve2}
  P' = - Q (t) - A^{T}(t) P - P A (t) + P S(t) P  , \qquad
  P(T)=Q_{T},
\end{equation}
if the time-dependent matrices $Q (t), A (t), S(t)$ are
continuous,  $Q (t), S(t)$ are symmetric and non negative for
$t\in[0,T]$ and $Q_T$ is a symmetric non negative matrix.

This is a very important property because the matrix $P(t)$ is
coupled with the equation for the state vector and plays an
important role on its stability. Standard methods do not preserve
the positivity property of the matrix $P$. However, some of the
splitting methods we have presented in this work will preserve
this property.

Notice that in the split (\ref{eq:Separable2tb}) one freezes the
time in the matrix RDE and solves exactly the corresponding
autonomous equation
\begin{equation}\label{eq:Preserve3}
  P' = - Q (t_1) - A^{T}(t_1) P - P A (t_1) + P S(t_1) P  , \qquad
  P(t^*)=P^*,
\end{equation}
where $Q(t_1),S(t_1)$ are symmetric non negative matrices. Then if
$P^*$ is a symmetric non negative matrix the numerical solution
$P_{ap}(t)$ is also a symmetric matrix and $P_{ap}(t)\geq 0, \
t\in[t^*-h,t^*]$.

Unfortunately, it is well known that splitting methods of order
higher than two necessarily have at least one coefficient $a_i$
negative and at least one coefficient $b_i$ negative, and then
splitting methods of order greater than two can not guarantee
positivity (a similar result was obtained for the preservation of
positivity of commutator-free methods studied in
\cite{blanes12mif}). One expects, however, a better preservation
of this property with respect to standard methods. In a splitting
method we can choose either the coefficients $a_i$ or $b_i$ to
advance the RDE. If positivity is an important property, it is
convenient to advance the RDE with maps associated to the
coefficients having the smaller negative value. For example, we
have interchanged the coefficients $a_i$ and $b_i$ in the 6-stage
fourth-order method (\ref{4thOrder}) as given in
\cite{blanes02psp} so $\ \min\{b_i\}>\min\{a_i\} \ $ because the
RDE is advanced with the coefficients $b_i$ in the algorithm.

\section{Differential games with $N$ players} \label{section2.1}
 \label{DifferentialGames}


Let us now consider the problem of differential games  with $N$
players given by the equations
\begin{equation}\label{eq:GamesNx}
 x'(t) =  A (t) x(t) + \sum_{i=1}^{N} B_{i} (t) u_{i}(t),  \qquad x(0) =
 x_0 \, ,
\end{equation}
and the  quadratic cost function
\begin{eqnarray}\label{eq:GamesNJ}
 & & J_{i}= x^{T}(T) \, Q_{iT} \, x(T)  \\
 & &  + \int_{0}^{T} \left\{ x^{T}(t) Q_{i} (t) x(t) +
u_{i}^{T}(t) R_{ii} (t) u_{i}(t)+ \sum_{j=1,j \neq i}^{N}
u_{j}^{T}(t) R_{ij} (t) u_{j}(t) \right\} dt  \nonumber
\end{eqnarray}
$ i=1, \ldots,N$.
From \cite{SH69}, for a zero sum game it is necessary that $
R_{ij} \neq 0 $, $ i \neq j $, but in a non-zero sum game it is
natural to choose $ R_{ij} = 0 $, $ i \neq j $, because in the
most frequent applications the cost function of each player does
not contain the other player control. In this way, the quadratic
cost function $J_{i}$ depends only on the control $u_{i}$.

Let us first consider a non-cooperative non-zero sum game where
each player, in order to minimize their cost function, determines
his action in an independent way knowing only the initial state of
the game and the model structure.
Under these conditions the optimal controls are given by
\begin{displaymath}
u_{i} (t) \, = \, - \, R_{ii}^{-1} (t) \, B_{i}^{T} (t)  \, P_{i}
(t) \,  x(t) \, , \qquad i=1, \ldots,N,
\end{displaymath}
where the matrices $P_i$, satisfy the coupled matrix RDEs
\begin{equation}\label{N}
\begin{array}{c}
 P'_{i}  =   - Q_{i} (t) - A^{T} (t) P_{i} - P_{i} A (t)
  + \sum_{j=1}^{N} P_{i}  S_{j} (t)  P_{j} \,
  \end{array}
\end{equation}
with $S_{i} (t) = B_{i} (t) \, R_{ii}^{-1}  (t) \, B_{i}^T  (t)$,
$i = 1, \ldots,N$. If we denote by
\begin{equation}\label{1}
W(t)=   \left[ \begin{array}{c} P_{1} (t) \\ \vdots \\ P_{N} (t)
\end{array} \right]  \in  \mathbb{R}^{Nn \times n}
\, , \qquad C(t) = \left[ \begin{array}{c} -Q_{1} (t) \\ \vdots \\
-Q_{N} (t)
\end{array} \right] \in \mathbb{R}^{Nn \times N}\, ,
\end{equation}
\bigskip
\begin{equation}\label{2}
B(t) = \left[ \, -S_{1} (t) \ \ \cdots \ \ -S_{N} (t) \, \right]
\in \mathbb{R}^{n \times Nn} \, ,
\end{equation}
\bigskip
\begin{equation}\label{3}
D(t) \, = \, \left[
\begin{array}{cccc} A^{T}(t) & 0 & \cdots & 0
\\
0 & A^{T}(t) & \cdots & 0 \\
\vdots & \vdots & \ddots & \vdots \\
0 & 0 & \cdots & A^{T}(t)
\end{array} \right] \in \mathbb{R}^{Nn \times Nn}\, ,
\end{equation}
then, the coupled system (\ref{N}) can be written as
\begin{displaymath}
W'(t)=C(t)-D(t)W(t)-W(t)A(t)-W(t)B(t)W(t) \, .
\end{displaymath}
From \cite{JoPo,reid}, if we consider $ \, y(t) \in
\mathbb{R}^{(N+1)n \times n} \, $ the solution of the linear
equation
\begin{equation}\label{4}
y'(t) \, = \, \left[ \begin{array}{cr} A(t) & B(t) \\ C(t) & -D(t)
\end{array} \right] \, y(t) \,
 ; \quad \quad
  y(T)=\left[ \begin{array}{c} U(T) \\ V_{1} (T) \\ \vdots \\ V_{N} (T) \end{array} \right]
  = \left[ \begin{array}{c} I \\ Q_{1T} \\ \vdots  \\ Q_{NT}  \end{array} \right] \,
  ,
\end{equation}
where $Z=[V_1^T,\ldots,V_N^T]^T$ and
$Z_T=[Q_{1T}^T,\ldots,Q_{NT}^T]^T$,and
we have that
$P_i(t) \, = \, V_i(t) \, U(t)^{-1}.$

As in the previous case, we can integrate backward in time the
coupled RDEs (\ref{4}) with a highly accurate method. Once the
initial conditions are obtained, we can integrate forward in time
the whole systems, which provides us the numerical values of the
controls along the time. The system to be solved is:

\begin{eqnarray}
 \frac{d}{dt} \left[
 \begin{array}{c} U \\ Z
 \end{array} \right] & = &
 \left[
 \begin{array}{ccc} A(t) & B(t) \\
 - C_{} (t)& - D (t)  \end{array} \right]
 \left[
 \begin{array}{c} U \\ Z
 \end{array} \right], \quad  \quad
 \left[ \begin{array}{c} U(t_0) \\ Z (t_0) \end{array} \right]
  = \left[ \begin{array}{c} U_0 \\ Z_{0} \end{array} \right]
    \label{eq:RiccatiForward}  \\
  \frac{dx}{dt}  & = &  \Big( A(t) - \sum_{i=1}^{N} S_i(t) \, V_i(t) \, U(t)^{-1} \,\Big) x ,
  \qquad \ \ x(0) = x_0 \,   \label{eq:StateForwardN}  \\
 u_{i} (t) & = & - \, R_{ii}^{-1} (t) \, B_{i}^{T} (t)
   \, V_i(t) \, U(t)^{-1}\, x(t) \, .
\end{eqnarray}

This problem has formally the same structure as the LQ optimal
control problem, and the numerical integration can be carried out
using the same methods as previously. The main difference remains
on the size of the matrices, and then its computational cost. For
a large number of players it can be advantageous to consider
approximations for the exponentials associated to the coupled
RDEs.

\subsection{Integrators for the zero sum game}

From \cite{SH69}, for a zero sum game it is necessary that $R_{ij}
\neq 0 $, $ i \neq j $. In this way, the quadratic cost functions
$J_{i}$ depend on all controls $u_{i}$, $ i= 1, \ldots, N \, $.
For simplicity in the presentation, we consider the case of two
players. The coupled Riccati differential equations to be solved
are
\begin{displaymath}
\begin{array}{rcl}
P'_{1} & = &  - Q_{1}(t) - A^{T}(t) P_{1} - P_{1} A(t)
    + P_{1}S_{1}(t) P_{1} + P_{1} S_{2}(t) P_{2} + P_{2} S_{22}(t) P_{2} \,
   ,  \\
P'_{2} & = & -Q_{2}(t) - A^{T}(t) P_{2} - P_{2} A(t)
   + P_{2}S_{2}(t) P_{2} + P_{2} S_{1}(t)  P_{1} + P_{1} S_{11}(t) P_{1} \, ,
\end{array}
\end{displaymath}
with the final conditions, $  P_{1}(T)=Q_{1T} \, , \ \
P_{2}(T)=Q_{2T} \, ,$ wherein
\[
S_{ij} (t) \, = \, B_{i} (t) \, R_{ij}^{-1}  (t) \, B_{i}^T  (t)
\, ; \ \ i = 1, \, 2 \, .
\]
This problem can not be reformulated as a linear problem. This
means we can not use exponential methods as previously. To
integrate backward the coupled RDEs we can use a highly accurate
method to avoid the loss in the preservation of the qualitative
properties. A high order extrapolation method based on a symmetric
second order integrator, say $\Phi_h^{[b]}$, would be a good
choice. Then, this symmetric second order method can be used as
the map in the composition (\ref{eq:S2}) for the forward time
integration, and we can use the composition methods given in
(\ref{eq:ProdS2}).

\section{Numerical examples}
 \label{NumericalExamples}


In order to test the performance of the numerical methods
presented in this work we consider  the problem of air pollutant
emissions studied in \cite{kata} and we have generalized this
problem to the case of $N$ regions ($N$ players) and the constant
parameters are replaced by time dependent functions
\begin{equation}\label{co2}
   x'= - a(t) x + b(t) \sum_{i=1}^N u_{i}(t)  ; \qquad x(0)=x_{0}.
\end{equation}
Here $x(t)\in\mathbb{R}$ is the excess of the pollutant in the
atmosphere, $u_{i}(t)\in\mathbb{R}$, $i=1, \ldots,N$, denote the
emissions of each region and $a(t)$, $b(t)$ are positive functions
related with the intervention of the nature on environment. The
cost functions to minimize is given by
\begin{equation}\label{costex}
J_{i} = 
 \int_{0}^{T} e^{- \rho t} \left\{ c_{i}(t)
u_{i}^{2} (t) + d_{i}(t) x^{2} (t) \right\} \, dt \, , \ \ i = 1 ,
\ldots,N \, ,
\end{equation}
where $c_{i}(t)$, $d_{i}(t)$, are positive functions  related to
the costs of emission and pollution withstand respectively, and $
\rho $ is a refresh rate. With an appropriate change of constants,
the problem (\ref{co2})-(\ref{costex}) can be applied to other
situations, such as financial problems.

In our notation, $R_{ij}=0$, $i \neq j $, $Q_{iT}=0$, $Q_{i} (t) =
d_{i}(t) e^{-\rho t}$, $R_{ii} = c_{i}(t) e^{-\rho t}$, and $
S_{i} = b(t)^{2} e^{\rho t}/c_{i}(t)$,
$ i=1, \ldots,N$.


Then, 
we have:
\begin{itemize}
\item  the equation (\ref{4}) with the data
\begin{displaymath}
 y'(t) = K(t) y(t) \, , \ \qquad
 y(T) = \left[ 1 , \, 0, \, \ldots \, , 0 \right]^T \, ;
\end{displaymath}
\begin{equation}\label{Ex2Ric}
y(t) = \left[ \begin{array}{c} u(t) \\ v_{1} (t) \\ \vdots \\
v_{N} (t) \end{array} \right] \, ; \ \ K(t) = \left[
\begin{array}{cccc}
- a(t) & - \displaystyle \frac{b(t)^{2}}{c_{1}(t)} e^{\rho t} &
\cdots &
- \displaystyle \frac{b(t)^{2}}{c_{N}(t)} e^{\rho t} \\
- d_{1}(t) e^{-\rho t} &  a(t) & \cdots & 0 \\
 \vdots &  \vdots & \ddots & \vdots \\
- d_{N}(t) e^{-\rho t} & 0 & \cdots &  a(t)
\end{array} \right] \, ,
\end{equation}
that must be solved from $t=T$ to $t=0$.
\item To integrate forward in time the IVP
\begin{eqnarray}\label{PhiEx2}
 y'(t) & = & K(t) y(t) \, , \ \qquad
 y(0) = y_0 \, , \\
 x' & = & \left( - a(t) - b(t)^{2} e^{\rho t} \sum_{i=1}^N
\frac{1}{c_{i}(t)} \, \frac{v_{i} (t)}{u (t)}  \right) x \, , \ \
x(0)=x_0 \, , \\
 u_{i} (t) & = & - \, \frac{1}{c_{i}(t)} \, e^{\rho
t} \, b(t) \, \frac{v_{i} (t)}{u (t)} \,x(t) \, , \qquad \ i  = 1,
\ldots,N.
\end{eqnarray}
\end{itemize}

In our numerical experiments we take $T=1$, initial conditions,
$x_0=10$ and consider the case of $N=10$ players, and we will take
different choices for the functions and parameters in order to
analyze the performance of the methods on different conditions.

The new methods, denoted by SP$n$, where $n$ is the order of the
method for the methods given in Table~\ref{tab.1}, will be tested
versus the following standard numerical methods:
\begin{itemize}
\item RK4: The well known 4-stage fourth-order Runge-Kutta method.
\item ODE45: The variable step and variable order algorithm \verb"ode45"
  implemented in \verb"Matlab".
\item  ODE113: The variable step and variable order algorithm \verb"ode113"
  implemented in \verb"Matlab".
\end{itemize}

We first consider the different cases for the autonomous problem
and next we consider other cases in which different explicitly
time-dependent functions are taken into account.

\paragraph{The autonomous problem}

We first consider the case where the RDE is autonomous, which
corresponds to the case $\rho=0$ and the functions $a,b,c_i,d_i$
take constant values, so $K$ in (\ref{Ex2Ric}) is a constant
matrix and the solution to get the initial conditions $y(0)$ is
given by
\[
  y(0) = e^{-TK}y(T).
\]
The initial conditions for the matrix RDE is computed up to round
off accuracy and the performance of the methods is measured by
taking into account the forward time integration. The solution at
$t=T$ satisfies that $u_i(T)=0, \ i=1,\ldots,N$, and this is
exactly satisfied, up to round off error, by the splitting methods
so we measure the error as $|x(T)-x_{ap}(T)|$, where $x_{ap}$ is
the approximated solution obtained by the numerical integrators.

\begin{enumerate}
\item We first take the following values for the parameters:
\[
 a=b=1, \qquad c_i=\frac1{d_i} = \frac{10+i}2.
\]
We measure the accuracy versus
  the number of function evaluations for each method when the numerical
  integration is carried out using different values of the time step
  (or different values for the absolute and relative tolerances for the methods ODE45 and
  ODE113. In the numerical experiments we consider
  \verb"AbsTol"=$10^{-i}$, \verb"RelTol"=$10^{1-i}, \ i=0,1,2,\ldots$
  For the splitting methods we measure the cost as the
  number of stages. Obviously, the splitting methods require the evaluation
  of several exponentials of matrices and the
  computational cost will depend on the problem. However, for the autonomous
  case most exponentials can be computed at the beginning and be stored to be used
  repeatedly along the integration.


\begin{figure}
  \includegraphics[width=0.75\textwidth]{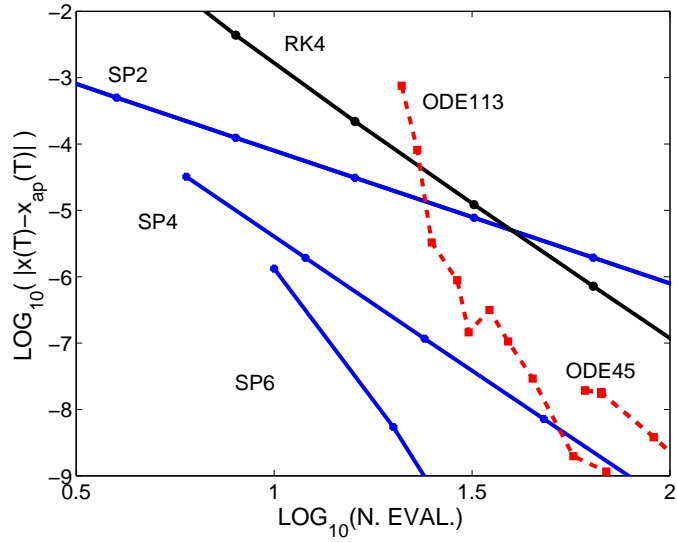}
\caption{ Error  $|x(T)-x_{ap}(T)|$ versus the number of
evaluations for the autonomous case with 10 players and $a=b=1, \
c_i=1/d_i=(10+i)/2$.}
\label{fig1}       
\end{figure}

  Figure~\ref{fig1} shows the
  results obtained. We observe that the splitting methods show
  more accurate results at the same number of evaluations. The performance of the new splitting methods could
  also be improved if the algorithms were implemented with variable order and
  variable step, as it is the case of ODE45 and ODE113.
\item We repeat the same experiment but for the case where
\[
  a=2, \qquad b=1, \qquad c_i=\frac1{d_i}=\frac{100+i}2.
\]
 This problem can be considered as a
  perturbed problem and we study the performance of the splitting
  methods (4,2) and (8,4) for the split given by
  (\ref{eq:PertForwardP}) and (\ref{eq:RiccatiForwardP2}),
  which are tailored for this class of problems.

We solve separately the dominant part of the system, given by the
equations
\begin{eqnarray}\label{}
  \left[ \begin{array}{c} u \\ v_{1} \\ \vdots \\
 v_{N} \end{array} \right]'
 & = & \left[    \begin{array}{cccc}
 - 2 &  0 & \cdots & 0 \\
   0 &  2 & \cdots & 0 \\
 \vdots &  \vdots & \ddots & \vdots \\
 0 & 0 & \cdots &  2
  \end{array} \right]
  \left[ \begin{array}{c} u \\ v_{1}  \\ \vdots \\
      v_{N} \end{array} \right], \qquad
  \left[ \begin{array}{c} u(t_n) \\ v_{1} (t_n) \\ \vdots \\
      v_{N} (t_n) \end{array} \right]=
  \left[ \begin{array}{c} u_n \\ v_{1,n}  \\ \vdots \\
      v_{N,n} \end{array} \right]
  \nonumber  \\
 x' & = & \left( - 2 - e^{\rho t} \sum_{i=1}^N
\frac{1}{c_{i}} \, \frac{v_{i} (t)}{u (t)}  \right) x \, , \qquad
x(0)=x_0 \, ,
\end{eqnarray}
where the RDE has trivial solution, which we plug into the
equation of the state vector
\begin{eqnarray}\label{}
  u(t) & = & e^{-2(t-t_n)}u_n  \\
  v_i(t) & = & e^{2(t-t_n)}v_{i,n}, \qquad i=1,\ldots,n  \\
 x' & = & \left( - 2 - e^{\rho t} e^{4 (t-t_n)} \sum_{i=1}^N
\frac{1}{c_{i}} \, \frac{v_{i,n}}{u_n}  \right) x \, .
\end{eqnarray}
The equation for the state vector has exact solution since it is a
scalar equation. In case it was a matrix equation, one could use,
for example, a Magnus integrator, where the exponentials from the
RDE can be computed for one time interval, and then used for all
the integration. Figure~\ref{fig2} shows the  results obtained.
The schemes SP4 and SP6 lead to slightly worst results and are not
showed.
\end{enumerate}


\begin{figure}
  \includegraphics[width=0.75\textwidth]{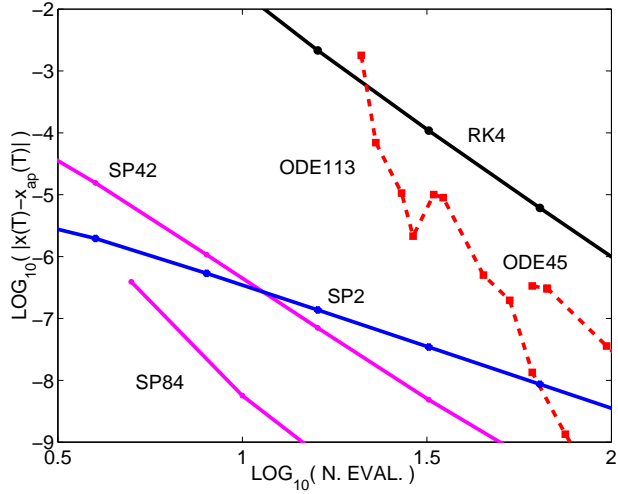}
\caption{ Error  $|x(T)-x_{ap}(T)|$ versus the number of
evaluations for the autonomous case with 10 players and $a=2, b=1,
\ c_i=1/d_i=(100+i)/2$.}
\label{fig2}       
\end{figure}

\paragraph{The non autonomous problem}

We first consider the one player case, $N=1$, with the following
choice for the functions and parameters
\[
 a(t)=2+\tanh\left(5(t-\frac12)\right), \quad \rho=\frac1{10},
 \quad
 b=1, \quad c_1=\frac1{d_1}=\frac{11}2.
\]
where $a(t)\in[1,3]$, and we repeat the same experiments replacing
the values of $c_1,d_1$ by
\begin{equation}\label{eq:coefsEj2b}
  c_1=\frac1{d_1}=\frac{101}2
\end{equation}
which makes the system closer to a near integrable systems, in
order to study the performance of the splitting methods in this
case.  We consider the split (\ref{eq:Separable2tb}) and the
algorithm (\ref{eq:algorithm1}). The case of one player
corresponds to a LQ optimal control problem and, as already
mentioned the solution of the RDE has to be a positive function.
To show the superiority of the splitting methods when this
property is important to be preserved, we measure the value of
$P(T)=\frac{v(T)}{u(T)}$ for each method and choice of the time
step. If this value is negative and smaller than a tolerance
value, which we take $-10^{-8}$ (i.e. if $P_{ap}<-10^{-8}$) we
marked this result with a circle.

\begin{figure}
  \includegraphics[width=0.5\textwidth]{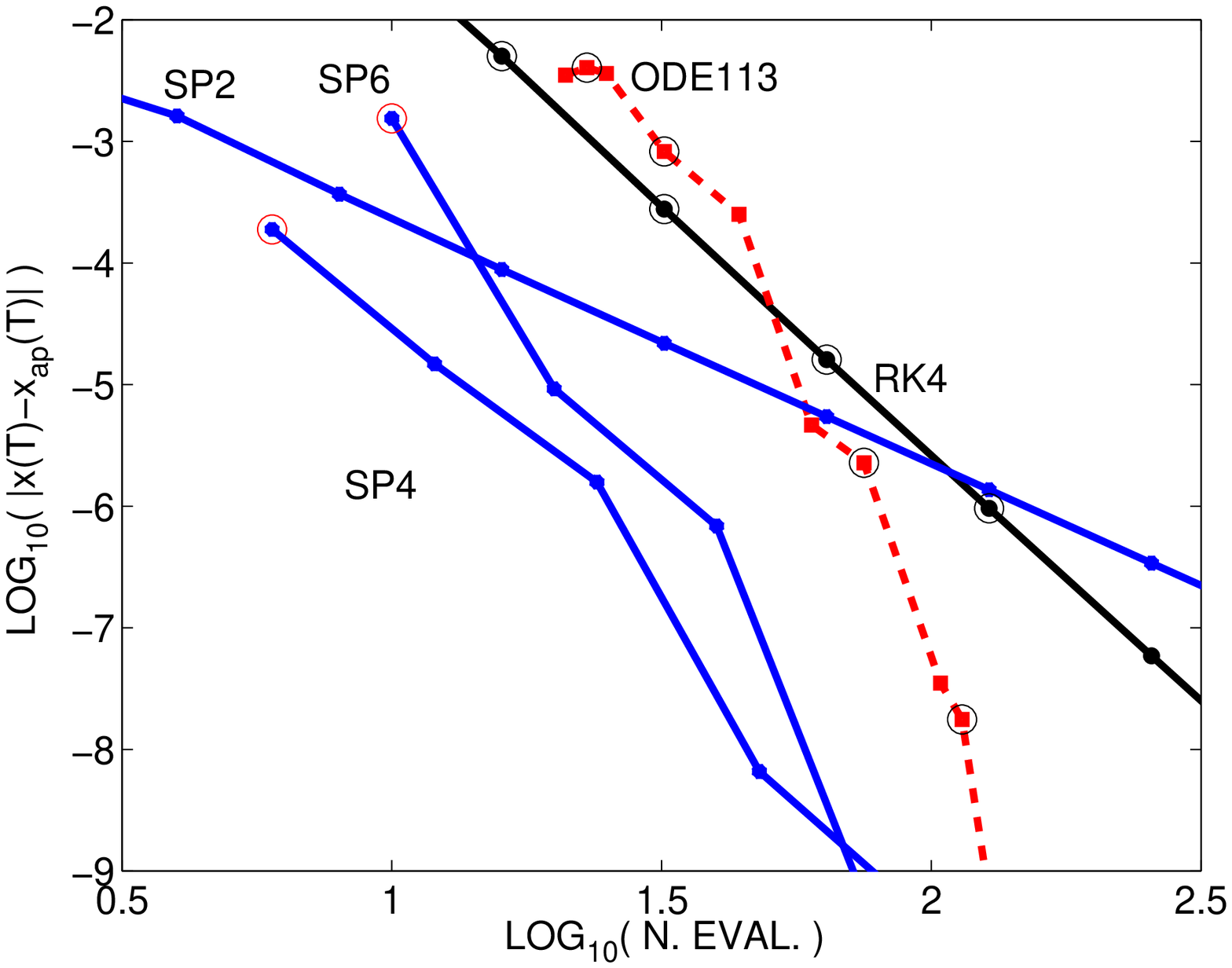}
  \includegraphics[width=0.5\textwidth]{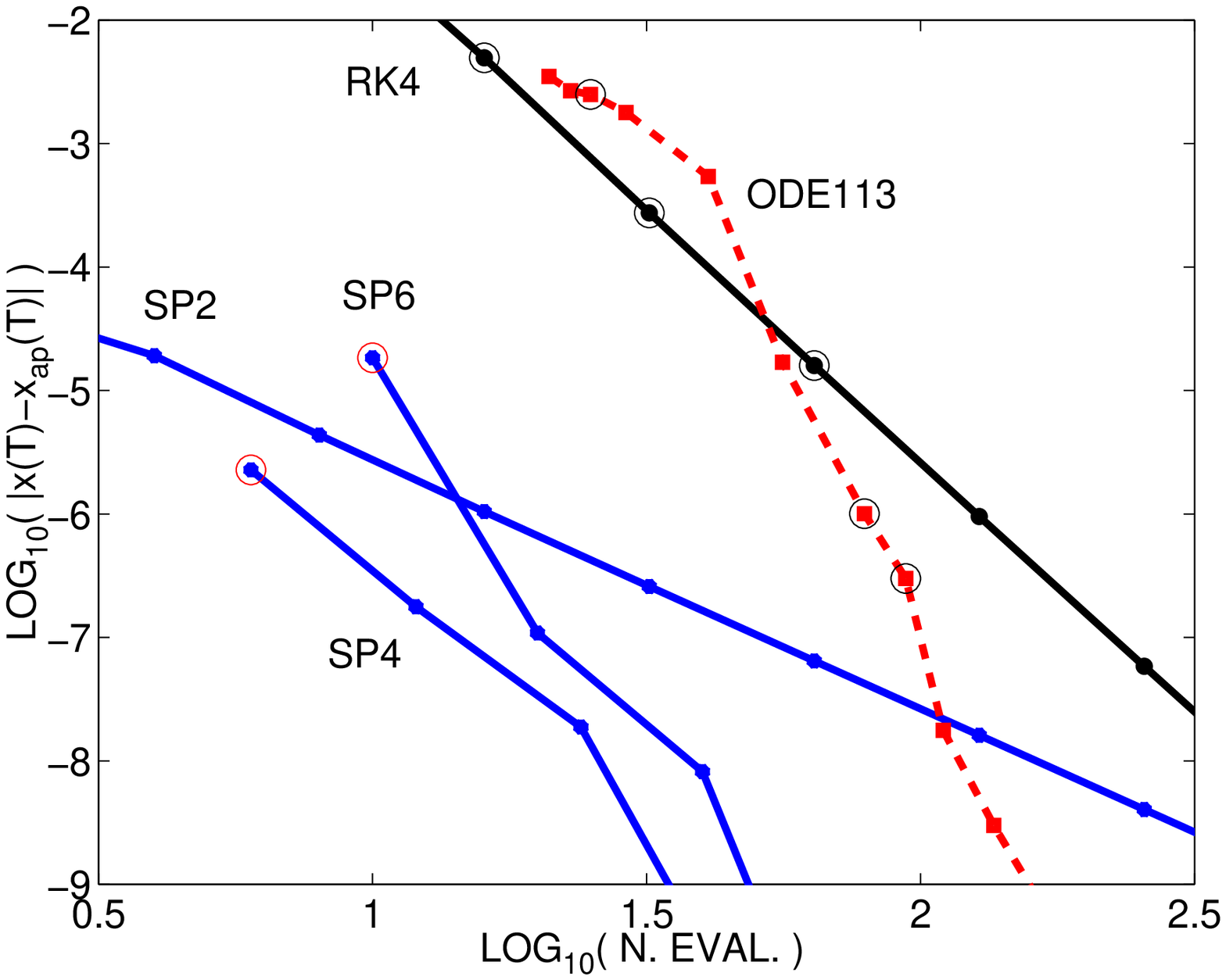}
\caption{ Error  $|x(T)-x_{ap}(T)|$ versus the number of
evaluations for the autonomous case with 1 player and
$a(t)=2+\tanh(5(t-1/2)), \ b=1, \ \rho=1/10$ and the following
values for $c_1,d_1$: Left panel  $c_1=1/d_1=11/2$; right panel
$c_1=1/d_1=101/2$. Solutions where $P_{ap}(T)< -10^{-8}\ $ are
marked with a circle.}
\label{fig3}       
\end{figure}

Figure \ref{fig3} illustrates the results obtained. The second
order splitting methods preserves positivity for all time steps
and the fourth- and sixth-order splitting methods do not preserve
this property only for the largest time step, contrarily to the
standard RK method or the methods implemented in \verb"Matlab". We
observe that when the off diagonal coefficients of the RDE are
small, the superiority of the splitting methods is even higher. We
have repeated the numerical experiments with a higher number of
players, and similar results are obtained.


%


\section{Conclusions}
 \label{Conclusions}

We have considered the numerical integration of linear-quadratic
optimal control problems and $N$-player differential games.
 These problems require to solve, backward in
time, non-autonomous matrix Riccati differential equations which
are coupled with the linear differential equations for the dynamic
state. The solution of both problems allows to obtain the optimal
control to be applied on the system in order to get the desired
target. The system of equations which describe the problems have a
particular algebraic structure which makes the solution to have
some qualitative properties. We present high order explicit
geometric integrators which preserve most qualitative properties.
In particular, we analyze the positivity of the solution for the
associated matrix Riccati differential equation, and we observe
that the methods presented show a better preservation for this
property. The coupled system of equations correspond to a
non-linear boundary value problems which usually are solved by
implicit and computationally very expensive algorithms. The
methods proposed are explicit and are not iterative methods.

The methods proposed are high-order explicit geometric integrators
which consider the time as two new coordinates. This allows us to
integrate the whole system forward in time while preserving the
most relevant qualitative structure of the exact solution. This
allows an immediate evaluation of the control forward in time and
then available for real time integrations. The numerical examples
considered showed the performance of the proposed methods as well
as the good preservation of some of the qualitative properties.
Similar ideas could be used for solving non-linear optimal control
problems after linearization \cite{blanes12mif}, and using an
iterative process, being this problem under consideration at this
moment.

\section*{Acknowledgement}


The  author wishes to thank the University of California San Diego
for its hospitality where part of this work was done. He also
acknowledges the support of the Ministerio de Ciencia e
Innovaci\'on (Spain) under the coordinated project
MTM2010-18246-C03.


\end{document}